\documentclass[11pt,a4paper]{article}

\title{Convex polytopes and linear algebra}

\author{Jonathan Fine\relax
\thanks{203 Coldhams Lane, Cambridge, CB1 3HY, England.  
\quad E-mail: \texttt{j.fine@pmms.cam.ac.uk}}
}
\date{1 October 1997}

\newcommand\PDelta{{\bf P}_\Delta}
\newcommand\bfR{{\bf R}}

\newcommand\bound{{\rm d}}
\newcommand\im{\mathop{\rm im}\nolimits}
\newcommand\ord{\mathop{\rm ord}\nolimits}
\newcommand\codim{\mathop{\rm codim}\nolimits}
\newcommand\bardiagram[1]{\leavevmode\hbox{\texttt{#1}}}
\newcommand\bibrule{\rule{2pc}{0.4pt}}

\begin{document}

\maketitle

\begin{abstract}\noindent
This paper defines, for each convex polytope $\Delta$, a family
$H_w\Delta$ of vector spaces.  The definition uses a combination of linear
algebra and combinatorics.  When what is called exact calculation holds,
the dimension $h_w\Delta$ of $H_w\Delta$ is a linear function of the flag
vector $f\Delta$.  It is expected that the $H_w\Delta$ are examples, for
toric varieties, of the new topological invariants introduced by the
author in \emph{Local-global intersection homology} (preprint
alg-geom/9709011).
\end{abstract}

\section{Introduction}

The goal, towards which this paper is directed, is as follows.  Suppose
$\Delta$ is a convex polytope.  One wishes to construct from $\Delta$
vector spaces whose dimension is a combinatorial invariant of $\Delta$. 
The smaller the dimension of these spaces, the better.  The convex
polytope $\Delta$ has both a linear structure, due to the ambient affine
linear space, and a combinatorial structure, due to the incidence
relations among the faces.  The construction in the paper uses both
structures to produce many `vector-weighted inclusion-exclusion formulae',
to each one of which corresponds a complex of vector spaces.  When such a
complex exactly computes its homology (a concept to be explained later)
the result is a vector space of the type that is sought.  The proof of
exact calculation, which is not attempted in this paper, is expected to be
difficult.

In a special case, part of this problem has already been solved.  If
$\Delta$ has rational vertices then from $\Delta$ a projective algebraic
variety $\PDelta$ can be constructed, and the middle perversity
intersection homology (mpih) Betti numbers $h\Delta$ are combinatorial
invariants of $\Delta$, by virtue of the Bernstein-Khovanskii-MacPherson
formula \cite{bib.JD-FL.IHNP,bib.KF.IHTV,bib.RS.GHV}, that express the
$h_i$ as linear functions of the flag vector $f\Delta$ of $\Delta$.  In
this case Braden and MacPherson (personal communication) have proved an
exact calculation result.  Their proof relies on deep results in algebraic
geometry, and in particular on Deligne's proof \cite{bib.PD.WC1,bib.PD.WC2}
of the Weil conjectures.  Elsewhere \cite{bib.JF.LGIH}, the author has defined
local-global intersection homology groups.  The construction in this paper
corresponds to the extension and unwinding (\cite{bib.JF.LGIH}, formulae
(8)--(10)) of the extended $h$-vector defined in that paper.  This
correspondence, which is an exercise in combinatorics, is left to the
reader.  It might also be presented elsewhere. This paper has been written
to be independent of \cite{bib.JF.LGIH}.

Central to this paper is the study of flags.  They too have linear and
combinatorial strucutre.  A \emph{flag} $\delta$ on a convex polytope
$\Delta$ is a sequence
\[
    \delta = ( \delta_1 \subset \delta_2 \subset \dots
                \subset \delta_r \subset \Delta )
\]
of faces $\delta_i$ of $\Delta$, each strictly contained in the next.  The
\emph{dimension} $d = \dim \delta$ is the sequence
\[
    d = ( d_1 < d_2 < \dots < d_r < n = \dim \Delta )
\]
of the dimensions $d_i$ of the \emph{terms} $\delta_i$ of $\delta$. 
Altogether there are $2^n$ possible flag dimensions.  The number $r$ is
the \emph{order} $\ord \delta$ of the flag $\delta$ (and of the dimension
vector $d$).  If $\delta$ is a face of $\Delta$, the order one flag whose
only term is $\delta$, namely $(\delta\subset\Delta)$, will also be
denoted by $\delta$.  The empty flag will be denoted by $\Delta$.

Similarly, if $V$ is a vector space then a flag $U$ on $V$ is a sequence
\[
    U = ( U_1 \subset U_2 \subset \dots \subset U_r \subset V )
\]
of subspaces $U_i$ of $V$, each stricly contained in the next.  The
\emph{dimension} $d=\dim V$ is the sequence
\[
    d = ( d_1 < d_2 < \dots < d_r < n = \dim V )
\]
where now the $d_i$ are the dimensions of the $U_i$.  Throughout $V$ will
be the vector space $\langle\Delta\rangle$ spanned by the vectors that lie
on $\Delta$.  Each face $\delta_i$ similarly determines a subspace
$U_i=\langle\delta_i\rangle$ of $V$.  Thus, each flag $\delta$ of faces on
$\Delta$ determines a flag $\langle\delta\rangle$ of subspaces of
$V=\langle\Delta\rangle$.

The construction of this paper is, in general terms, as follows.  Suppose
$U$ is a flag on $V$.  From $U$ many vector spaces can be constructed. 
This paper constructs for each $U$, and for each $w$ lying in an as yet
unspecified index set, a vector space $U_w$.  Now let $U'$  be obtained
from $U$ by deleting from $U$ one of its terms $U_i$.  It so happens that
this deletion operator induces a natural map
\[
    U_w \to U'_w
\]
between the associated $w$-spaces.  (In the simplest case, which
corresponds to mpih, all the $U_w$ are subspaces of a single space $V_w$,
and the maps are inclusions.  In general there is no such a global space
$V_w$.)  Now suppose $U'$ is obtained from $U$ by deleting two or more
terms.  By choosing an order for the deletion of these terms, a map 
$U_w \to U'_w$ can be obtained.  The single-deletion map is natural (or
geometric) in that the induced multiple deletion maps are independent of
the choice of deletion order.

Suppose now that such spaces $U_w$ have been defined for every flag on
$V$, and that $\Delta$ is a convex polytope whose vector space
$\langle\Delta\rangle$ is $V$.  Each flag $\delta$ on $\Delta$ thus
determines a flag $\langle\delta\rangle$ on $V$, and thus a vector space
$\langle\delta\rangle_w$, or $\delta_w$ for short.  These vector spaces
will be assembled into a complex, according to the order $r$ of $\delta$. 
Define the space $\Delta(w,r)$ of \emph{$w$-weighted $r$-flags} to be the
direct sum of the $\delta_w$, where $\delta$ has order $r$.  Each vector
in $\Delta(w,r)$ can be thought of as a formal sum $\sum v_\delta
[\delta]$, where $[\delta]$ is a formal object representing an $r$-term
flag $\delta$, and where the coefficient $v_\delta$ is drawn from the
vector space $\delta_w$.)  By the assumptions of the previous paragraph,
the deletion operator on flags determines a differential
\[
    \bound: \Delta(w,r) \to \Delta(w,r-1)
\]
and so induces a complex
\[
    0 \to \Delta(w,r) \to  \Delta(w,r-1) \to \dots
     \Delta(w,1) \to \Delta(w,0) \to 0
\]
of vector spaces.  (As is usual, $\bound=\sum (-1)^{i+1}\partial_i$, where
$\partial_i$ is the operator induced by deletion of the $i$-th term. 
Because the maps $U_w\to U'_w$ are natural, $\bound^2$ is zero, and thus
one indeed has a complex.)

\emph{Exact calculation} is when this complex is exact except at one
point, say $\Delta(w,j)$.  In that case the homology $H_w\Delta$ at that
point is a suitable alternating sum of the dimensions of the
$\Delta(w,i)$, and thus of the $\delta_w$.  By construction the dimension
of $\delta_w$ will depend only on the dimension vector $d$ of $\delta$,
and thus (provided exact calculation holds) one has that the dimension
$h_w\Delta$ of $H_w\Delta$ is a linear function of the flag vector
$f\Delta$ of $\Delta$, and so is a combinatorial invariant.  (This is
because the $f_d\Delta$ component of $f\Delta$ counts how many $d$-flags
there are on $\Delta$.  Each $f_d\Delta$ contributes, up to an alternating
sign, the quantity $\lambda_d=\dim\delta_w$ to $h_\Delta$, where
$\delta_w$ is a coefficient space due to any flag $\delta$ of dimension
$d$.)

Convex polytopes are not the only combinatorial objects for which the
concept of a flag can be defined.  In \cite{bib.JF.QTHGFV,bib.JF.SFV}, the
author defines flag vectors for $i$-graphs, or more generally any object
that is a union of cells (or edges), and which can be shelled.  There
seems to be no reason why the general form of the construction described
here cannot also be applied in this new context.  This is not to say the
the proper choice of the vector spaces $\delta_w$ to be associated to the
flags $\delta$ is not expected to be a deep question.  Both for convex
polytopes and for $i$-graphs there are subtle and presently unknown
combinatorial inequalities on the flag vectors.  An exact homology theory,
using vector spaces such as the $\delta_w$, is the only method the author
can envision, that will lead to the proof such inequalities.

The exposition is organised as follows.  The next section (\S2) give the
definition of exact homology.  The deletion operator applied to flags
yields the flag complex~(\S3).  Next comes a complex~(\S4) that
corresponds to middle perversity intersection homology.  The local-global
variant is more complicated, and is the substance of the paper.  First the
coefficient spaces to be used are defined~(\S5), then the maps between
them~(\S6), and finally the local-global homology~(\S7). 

The purpose of \S\S3--7 is to present the definition as the result of a
study of the geometric resources and constraints.  Conversely, in \S8
decorated bar diagrams are used to reformulate the definition, and allow
the properties to be demonstrated, in a concise manner.  Finally, \S9
provides a wider discussion of what has been done, and what remains to be
done.

The reader may at first find \S\S5--7 somewhat abstruse.  They contain a
study of the linear algebra of a segmented flag of vector spaces.  Once
the problem is understood, the key definitions come out in a fairly
natural way.  In \S8, the same definitions are presented, but this time
via coordinates.  Here, the definitions are clear, but may appear somewhat
arbitrary.  Each point of view informs the other.  This paper attempts to
show how the definitions arise naturally out of the logic of the
situation, and thus places \S\S5--7 before \S8.  The reader may wish to
reverse this order.

The basic devices are the combinatorics of flags, and linear algebra.  The
exterior algebra on a vector space is widely used.  The fibre of the
moment map from a projective toric variety to its defining polytope is
always a product of circles, and so the homology of the fibre is
isomorphic to an exterior algebra.  This fact is the beginning of the
connection between the linear algebra of flags and the existence of cycles
on the toric variety.  The reader does not need to know this.

Throughout this paper $\Delta$ will be a convex polytope of some fixed
dimension~$n$, and $V$ will be the vector space spanned by the vectors
lying on $\Delta$.  It will do no harm to think of $\Delta$ as lying in
$V$.  

\section{Exact Homology}

This section explains the concept of exact homology.  Suppose that a
sequence
\[
    A =
    ( 0 \to A_n \to A_{n-1} \to \dots \to A_1 \to A_0 \to 0 )
\]
of vector spaces is given, together with a map
\[
    \bound:A_i \to A_{i-1}
\]
between successive terms.  If, for each $i$, the composite map
\[
    \bound \circ \bound =  \bound^2 : A_i \to A_{i-2}
\]
is the zero map, then $A$ is called a \emph{complex}, and $d$ its
\emph{boundary map}.  If $A$ is a complex then the statement
\[
    \im (\bound:A_{i+1} \to A_i ) 
    \> \subseteq \>
    \ker (\bound:A_i \to A_{i-1})
\]
restates the condition $\bound^2=0$, and the quotient of the above kernel
by the image is called the $i$-th homology $H_iA$ of the complex $A$. This
formalism originated in the definition, via chains and cycles, of the
homology groups $H_iX$ of a topological space $X$.  There, most of the
homology groups were expected to be non-zero. The present use will be
different.

Suppose $A$ is a complex.  If all the homology groups $H_iA$ are zero
(i.e.~at each $A_i$ the kernel and image are equal) then $A$ is called a
(long) \emph{exact sequence}, or \emph{exact} for short. If $A$ is exact
then the alternating sum
\begin{equation}
\label{eqn:chi}
    \sum \nolimits _{i=0}^{n}\> (-1)^i \dim A_i
\end{equation}
of the dimensions of the $A_i$ (assumed finite) will be zero.  To prove
this, introduce in each $A_i$ a subspace $B_i$ that is a complement to
$\bound A_{i+1}$ in $A_i$.  The dimension of $A_i$ is the sum of the
dimensions of $\bound A_{i+1}$ and of $B_i$.  The exactness assumption
implies that the restricted form
\[
    \bound:B_i \to \bound A_i \subseteq A_{i-1}
\]
of the boundary map is an isomorphism.  Thus, the contributions of $B_i$
and $\bound A_i$ to the alternating sum are equal but opposite, and so the
result follows.

Suppose $A$ is a complex.  If $H_iA$ is zero then $A$ is said to be
\emph{exact at $A_i$}.  Now suppose that $A$ is known to be exact at all
its $A_i$ except perhaps one, say $A_r$.  In this case the alternating sum
(\ref{eqn:chi}) of the dimensions gives not zero, but the dimension of the
only non-zero homology $H_rA$ of the complex $A$, multiplied by $(-1)^r$. 
If a complex is exact at all but one location, we shall say that it
\emph{exactly computes} the homology at that location.

Now suppose that the $A_i$ are constructed from the convex polytope
$\Delta$, and that their dimensions depend only on the combinatorial
structure of $\Delta$.  The same will then be true for $H_rA$, provided
that for each $\Delta$ the complex exactly computes its homology at
$A_r$.  If this holds we shall say that $H_rA$ is an \emph{exact homology}
group of $\Delta$.  It is of course one thing to define a complex $A$, as
is done in this paper, and quite another to prove its exactness.  This is
expected to require new concepts and methods.

\section{The flag complex}

The convex polytope $\Delta$ has both a combinatorial structure (incidence
relations among faces) and a linear structure (the vectors lying on a face
$\delta$ span a subspace $\langle\delta\rangle$ of $V$).  In general both
structures will be used to define the complex $A$.  This section defines a
complex that uses the combinatorial structure alone.  The general case
will arise by allowing vectors to be used instead of numbers in the
construction that follows.

Recall the definition, in \S1, of a flag $\delta$ on $\Delta$, its
dimension vector $d$, and its order $r$.  Now suppose $\delta$ is a flag on
$\Delta$, of order $r$.  By removing one or more of the terms $\delta_i$
from $\delta$, new flags can be obtained, of lower order.  Let
$\partial_j$ be the \emph{deletion operator} that removes from a flag
$\delta$ the $j$-th term.

The operators $\partial_j$ do not commute.  Removing say the $2$nd term
from a sequence, and then say the $4$th, gives the same final result as
does first removing the $5$th and then the $2$nd.  This is because
removing the $2$nd term will cause the subsequent items to move down one
place in the sequence.  The equation
\[
    \partial_j \partial_k = \partial_{k+1} \partial_j
    \qquad \mbox{for $j<k$}
\]
is an example of the \emph{commutation law} for these deletion operators.

Now let $A_i$ consist of all formal weighted sums $x$ of the form
\[
    x = \sum \nolimits _ { \ord \delta = i } x_\delta [\delta]
\]
where the coefficients $x_\delta$ are numbers and, as indicated, the sum
is over all flags $\delta$ on $\Delta$, which have $i$ terms.  Here
$[\delta]$ denotes $\delta$ considered as a formal object.  Where
confusion will not then result, $\delta$ will be written in its place.
Thus, a vector $x$ in $A_i$ is a formal sum of $i$-term flags, with numeric
coefficients $x_\delta$.

The deletion operation $\partial_j$ induces a map $A_i \to A_{i-1}$.  (It
is zero if $j$ is larger than $i$.)  Now use the formula
\[
    \bound = \partial_1 - \partial_2 + \partial_3 + \dots
        + (-1)^{i-1}\partial_i + \ldots
\]
to define a map $\bound:A_i \to A_{i-1}$.  It follows immediately from the
commutation laws, that $\bound^2=0$, and so $A$ is a complex.  It will be
called the \emph{flag complex} of $\Delta$.  It depends only on the
combinatorial structure of $\Delta$.

\section{Global coefficient spaces}

Instead of using numeric coefficients for the formal sums that constitute
the space $A_i$, one could instead write
\begin{equation}
\label{eqn.sum.vdelta.delta}
    v = \sum \nolimits _ { \ord \delta = i } v_\delta \delta
\end{equation}
where now $v_\delta$ is to lie in a vector space $\Lambda(\delta)$ that is
in some way associated to $\delta$.  Thus, $A_i$ is again to be formal
sums of order $i$ flags, where now the coefficients are to be vectors
rather than numbers.  This section will describe the simplest way of
constructing such \emph{coefficient spaces} $\Lambda(\delta)$.

Suppose $\delta'$ is obtained from $\delta$ by the deletion of one or more
terms.  For such a definition to produce a complex, there must also be
a natural map
\[
    \Lambda(\delta) \to \Lambda(\delta')
\]
between the corresponding coefficient spaces.  One way to do this, which
will be used in this section, is to have this map be an inclusion.  Thus,
all the coefficients will lie in the global coefficient space
$\Lambda(\Delta)$ associated to the empty flag.  Each flag $\delta$ will
then define a subspace $\Lambda(\delta)$ of $\Lambda(\Delta)$.

For this to work, one must have that $\Lambda(\delta)$ is a subspace of
$\Lambda(\delta')$.  One way to do this is to have each individual face
$\delta_j$ in $\delta$ define a condition (or set of conditions) on
$\Lambda(\Delta)$.  Now define $\Lambda(\delta)$ to be those vectors in
$\Lambda(\Delta)$ that satisfy the condition(s) due to the faces
$\delta_j$ in $\delta$.  Provided the conditions due to $\delta_j$ depend
only on the face $\delta_j$ (and not on its location in $\delta$, or
whatever), it will follow automatically that $\delta'$ will provide fewer
conditions, and so it will be certain that $\Lambda(\delta')$ will
contain $\Lambda(\delta)$ as a subspace.

To summarise this section so far, suppose that a vector space
$\Lambda(\Delta)$ is given, and for each face $\delta$ of $\Delta$ a
subspace $\Lambda(\delta)$ is given (here $\delta$ stands for the flag
which has $\delta$ as its only term).  From this a complex $A$ can be
constructed.  First define $\Lambda(\delta)$ to be the intersection of the
spaces $\Lambda(\delta_j)$ associated to the terms $\delta_j$ of
$\delta$.  Next define $A_i$ to be all formal sums
(\ref{eqn.sum.vdelta.delta}, where the coefficients $v_\delta$ are to lie
in $\Lambda(\delta)$.  Finally, the boundary map
\[
    \bound (v_\delta \delta ) = v_\delta [\partial_1\delta]
                              - v_\delta [\partial_2\delta]
                               + \cdots
\]
is defined just as before.
To complete such a definition, one must provide a vector space
$\Lambda(\Delta)$, and derive from each face $\delta\subset\Delta$ a
subspace $\Lambda(\delta)$ of $\Lambda(\Delta)$.

Recall that $V$ stands for the span of the vectors lying on $\Delta$.  Let
$V^*$ be the dual space of linear functions.  That such a linear function
$\alpha$ is constant on $\delta$ (i.e.~zero on the vectors lying on
$\delta$) describes a subspace $\delta^\perp$ of $V^*$.

Exterior algebra will now be used.  Fix a degree $r$, and let
$\Lambda(\Delta)$ be the $r$-fold exterior product $\Lambda^r=
\bigwedge^rV^*$ of arbitrary linear functions on $V$.  The face $\delta$
defines a filtration of $\Lambda(\Delta)$ in the following way.  For
each decomposition $r=s+t$ one can take the span of expressions of the
form
\[
    \alpha_1 \wedge \alpha_2 \wedge \dots \wedge \alpha_{s}
    \quad \wedge \quad
    \beta_1 \wedge \beta_2 \wedge \dots \wedge \beta_{t}
\]
where the $\alpha_i$ are to vanish on $\delta$.  No conditions are placed
on the $\beta_i$.  The result is of course a subspace of
$\Lambda(\Delta)=\bigwedge^rV^*$.

To conclude this definition it is enough, for each face $\delta$ of
$\Delta$, to choose one of these subspaces of $\Lambda(\Delta)$.  One
would like the resulting complex to produce an exact homology group, a
matter which is presently not well understood, and which involves concepts
that lie outside the scope of this paper.

The condition
\[
    s > \codim \delta - s
\]
is satisfied by some smallest value of $s$.  (The \emph{codimension}
$\codim\delta$ is defined as usual to be $\dim\Delta-\dim\delta$.)  Use
this value to define for each face $\delta$ the subspace
$\Lambda(\delta)$.  If the resulting $s$ is greater than $r$, then
$\Lambda(\delta)$ is taken to be zero.

This choice of spaces corresponds to middle perversity intersection
homology (for the associated toric variety $\PDelta$, if it exists). There
is little doubt that this gives the correct choice of subspaces, for
reasons that will be discussed in the final section.

\section{Local coefficient spaces}

To produce a complex $A$ one requires a coefficient space
$\Lambda(\delta)$ for each flag $\delta$ on $\Delta$, and natural maps
$\Lambda(\delta) \to \Lambda(\delta')$ whenever $\delta'$ is obtained
from $\delta$ by the deletion of one or more terms.  The previous section
assumed the maps were inclusions (and so all the $\Lambda(\delta)$ were
subspaces of $\Lambda(\Delta)$).  This section will relax this assumption,
to obtain the coefficient spaces that will later be used to define further
complexes.  In the next section, the boundary map will be defined.

Suppose $V_1$ is a subspace of $V$.  A basic construction of the previous
section was to use $V_1$ to define subspaces (in fact a filtration) of the
$r$-fold exterior produce $\bigwedge^rV^*$.  Such a construction is in
fact forced upon us, provided we assume that the deletion operator on
flags induces inclusion, and also that $\Lambda(\Delta)$ is
$\bigwedge^rV^*$.

Relaxing this assumption allows the following.  Given $V_1\subset V$ one
can form the vector spaces $V_1$ and $V/V_1$ and then form the tensor
product
\[
    \bigwedge \nolimits ^{r_1} V_1^* 
    \otimes
    \bigwedge \nolimits ^{r_2} V_1^\perp
\]
of the corresponding exterior products.  Here, $V_1^\perp$ consists of the
$\alpha$ in $V^*$ that vanish on $V_1$.  It is of course naturally
isomorphic to $(V/V_1)^*$.

In this way $V$ can be broken into two or more segments, within each of
which the construction of the previous section can be applied.  For
example, each subspace $U$ of $V_1$ will determine a filtration of the
first factor $\bigwedge ^{r_1} V_1^*$ above. Similarly, if $U$ lies beween
$V_1$ and $V$, a filtration of the second factor $\bigwedge ^{r_2}
V_1^\perp$ will arise.

All the coefficient spaces $\Lambda(\delta)$ will be obtained in this way.
Given a flag $\delta$ use some (perhaps all or none) of its terms
$\delta_i$ to segment $V$.  This gives a tensor product of exterior
algebras.  For each factor choose a component, i.e.~a degree.  The
resulting space is used in the same way as $\Lambda(\Delta)$ was, in the
previous section.  Note that this space depends on the flag, or more
exactly the subflag used to segment $V$.  There is no longer one global
space, in which all the coefficient vectors lie.

Now choose a term $\delta_j$ of $\delta$.  This can be used to filter
$\bigwedge^{r_{k+1}} (V_{k+1}/V_k)^*$,  where $\delta_k$ is the largest
segmenting face contained in $\delta_j$.  (If there is none such, set
$k=0$ and use $\bigwedge^{r_1}V_1^*$ instead.)  Now, as before, use the
condition
\[
    s_j > \codim \delta_j - s_j
\]
to select a term in the filtration.  Here, however,
\[
    \codim \delta_j = \dim \delta_{k+1} - \dim \delta_j
\]
is to be the codimension of $\delta_j$ not within $\Delta$ but within its
segment.

The above filtration is to be applied for all the terms in the flag,
including those used to segment.  Such terms can produce of course only
a trivial filtration.  However, when
\[
    r_k > ( \dim \delta_{k} - \dim \delta_{k-1} ) - r_k
\]
holds, the whole of $\bigwedge^{r_{k}} (V_{k}/V_{k-1})$ is to be used, as
the coefficient space for $[\delta]$.  If the condition fails, zero is the
only coefficient to be used with $[\delta]$.

This concludes the definition of the coefficient spaces $\Lambda(\delta)$. 
Note that to specify such a space the following is required, in addition
to $\delta$.  One must select some (or all or none) of the faces of
$\delta$, to be used for segmentation.  One must also specify a degree
$r_k$ for each segment.  A more explicit notation might be
\[
    \Lambda ( \delta , -r , s )
\]
where $r$ is the sequence $(r_1, r_2, \ldots )$ of degrees, and $s$ gives
the subflag $\delta_s$ of $\delta$ that is used to segment $\delta$.  The
minus sign in $-r$ is to distinguish this notation from
$\Lambda(\delta,r,s)$, which will be introduced later.  As already noted,
each degree must be greater than half the length of the corresponding
segment, for the coefficient space to be non-zero.

\section{The boundary map}

Suppose that $\delta$ is a flag of $\Delta$, and that some segmentation
$s$ of $\delta$ is chosen, which breaks $\delta$ (and $V$) into $l$
segments. Suppose also that a multi-degree $r=(r_1,\dots, r_l)$ is given. 
The construction of the previous section yields a coefficient space
$\Lambda(\delta,-r,s)$.  Now suppose that $\delta'$ is obtained from
$\delta$ be removing one of the terms $\delta_j$ from $\delta$.  Provided
$\delta_j$ was not used to segment $\delta$, the space
$\Lambda(\delta',-r,s)$ will contain $\Lambda(\delta,-r,s)$.

One could obtain a complex by not allowing deletion to occur only at the
faces $\delta_j$ used to segment $\delta$, or rather setting the result of
deleting such a face to be zero.  The resulting complex and its homology
will not however properly speaking be an invariant of the polytope
$\Delta$.  Rather, for each choice of a segmenting flag $\delta_s$ one
will have a complex, and any invariant of $\Delta$ so defined will simply
be the direct sum of these flag contributions.

The purpose of this section is to define a map from $\Lambda(\delta)$ to
$\Lambda(\delta')$, where $\delta'$ is obtained from $\delta$ by the
removal of a segmenting face $\delta_i$.  This will have the effect of
producing a single global complex, that will in general be indecomposable. 
It can be thought of as a gluing together of the local complexes of the
previous paragraph.

Here is an example.  Suppose one has subspaces $V_1 \subset V_2 \subset
V$, where $V_i$ is the span $\langle\delta_i\rangle$ of vectors lying on
the face $\delta_i$.  The segmentation of $V$ due to $V_1$ will be
compared to that due to $V_2$.  In both cases, the raw materials are
firstly linear functions $\alpha$ defined on $V_i$, and secondly linear
functions $\beta$ vanishing on $V_i$ (and defined on the whole of $V$).

Now suppose that $\alpha_2$ is defined on $V_2$.  Because $V_1\subset
V_2$, the linear function $\alpha_2$ can be restricted to give $\alpha_1$
defined on $V_1$.  This is straightforward.  Now suppose that $\beta_2$
vanishes on $V_2$ (and is defined on $V$).  Again because $V_1 \subset
V_2$, the linear function $\beta_2$ also vanishes on $V_1$.  Thus there is
a linear map (restriction of range $\otimes$ relaxation of condition)
\begin{equation}
\label{eqn:rrmap}
    \bigwedge \nolimits ^r V_2^* 
        \otimes \bigwedge \nolimits ^s V_2^\perp
    \to
    \bigwedge \nolimits ^r V_1^* 
        \otimes \bigwedge \nolimits ^s V_1^\perp
\end{equation}
between the basic spaces associated to the two segmentations.

This map has an interesting relation to the conditions used to define the
coefficient spaces $\Lambda(\delta)$.  When the condition
\[
    s > \codim V_i - s
\]
holds, the coefficient space due to the flag $(\delta_i)$ will be one of
the above tensor products.  When the condition fails, the coefficient
space is zero.

Because $V_1$ is a subspace of $V_2$, the condition for $(\delta_1)$ is
more onerous than that for $(\delta_2)$, and so the map (\ref{eqn:rrmap})
is going in the wrong direction, to map the one coefficient space to the
other.

Regarding conditions however, the situation is different.  The trick is to
think of the conditions that define the coefficient spaces
$\Lambda(\delta_i)$ to be themselves subspaces of an exterior algebra
(or more exactly a tensor product of such).  If $U$ is a vector space of
dimension $l+m$, then each subspace of $\bigwedge^l U$ determines a
subspace of $\bigwedge^m U$, and vice versa.  This is via the
nondegenerate pairing
\[
    \bigwedge \nolimits ^l U \otimes \bigwedge \nolimits ^m U
    \to \bigwedge \nolimits ^{l+m} U 
\]
provided by the exterior algebra.  By a slight abuse of language, this
will be called duality.  (The value space $\bigwedge^{l+m}U$ is has
dimension one, but is not naturally isomorphic to $\bfR$.)

It is now necessary to formulate the conditions using this new point of
view.  Suppose $U\subset V$ is a subspace, and $W\subseteq \bigwedge^r
V^*$ is spanned by
\[
    \alpha_1 \wedge \dots \wedge \alpha_ l
    \quad \wedge \quad
    \beta_1 \wedge \dots \wedge \beta_ m
\]
where the $\alpha_i$ are to vanish on $U$.  Consider the space
$W'\subseteq\bigwedge^{r'}V^*$ spanned by
\[
    \alpha_1 \wedge \dots \wedge \alpha_{l'}
    \quad \wedge \quad
    \beta_1 \wedge \dots \wedge \beta_{m'}
\]
where $r+r'=\dim V$, $l+l'= \codim U + 1 $, and as before the $\alpha_i$
vanish on $U$; this is the subspace of the complementary component of the
exterior algebra, determined mutually by $W$.

The $\beta_i$ are arbitrary, and so by a further application of duality,
they can be dropped from the definition of the condition space.  Thus,
take the span of 
\[
    \alpha_1 \wedge \dots \wedge \alpha_{l'}
\]
for $l+l'=\codim U + 1$, $\alpha_i$ vanishing on $U$, as the condition
space for $W$.

We now return to the change of segmentation map (\ref{eqn:rrmap}).  As
already noted, the conditions for $V_2$ are less onerous that those for
$V_1$.  The condition space for $V_2$ is either zero (no conditions) or
$\bigwedge^{\codim V_2} V_2^\perp$ (zero is the only solution to the
conditions), and similarly for $V_1$.  The map (\ref{eqn:rrmap}) will in
either case respect these conditions.  (To make sense of this statement,
one should think of the conditions as being an ideal in the exterior
algebra, and then the image under (\ref{eqn:rrmap}) of the one ideal is
contained in the other.)

The goal of this section is now in sight.  It is the conditions that are
respected by the natural map that is due to change of segmentation, not
the coefficient spaces defined by the conditions.  Interpret each
coefficient vector $v_\delta$ as a linear function, taking values in a
tensor product of top-degree exterior products, that vanishes on the
condition space.

Recall that the example $V_1\subset V_2\subset V$ gives rise to a natural
map
\[
    \bigwedge \nolimits ^r V_2^* 
        \otimes \bigwedge \nolimits ^s V_2^\perp
    \to
    \bigwedge \nolimits ^r V_1^* 
        \otimes \bigwedge \nolimits ^s V_1^\perp
\]
which respects conditions.  Now let $v_\delta$ be a coefficient vector on
$V_1$, interpreted as a linear function on the range of the above map.

Using the map, this linear function on the range can be pulled back to
give a linear function on the domain.  Because $v_\delta$ vanishes on the
$V_1$ condition space, the pull-back vanishes on the $V_2$ condition
space.  By duality, the pull-back linear function is associated to a
unique coefficient vector $v'_\delta$ for the $V_2$ segmentation.  This is
an example of a \emph{change of segmentation component} of the boundary
map.  Note that it takes a $V_1$-segmentation coefficient to a $V_2$ such. 
In other words, under boundary the segmenting term(s) may move rightwards,
or in other words, increase in dimension.

To finish, there are some details to be taken care of.  First, although
the one-dimensional value spaces for $V_1$ and $V_2$ are not the same,
they are naturally isomorphic.  This is good enough.  Secondly, because
the map (\ref{eqn:rrmap}) is natural, the argument that shows $d^2=0$
works just as before.  Thirdly, the above argument has been applied only
to the conditions due to the segmenting faces.  The reader may wish to
show that it works also for the conditions, as imposed in the previous
section.  (An alternative way of defining the map and checking the
statements made will be outlined, as part of the discussion of bar
diagrams.)

The final matter concerns the degree of the coefficient $v_\delta$.  The
map (\ref{eqn:rrmap}) preserves the degree of (tensor products of)
exterior powers.  It induces the map on coefficients via duality and
pull-back, and so as (tensor products of) exterior powers the coefficients
$v_\delta$ on $V_1$ and $v_\delta'$ on $V_2$ will have different degrees. 
However, they will by definition have the same \emph{co-degree}, by which
is meant the amount by which they fall short of being of top degree.  For
this reason, in the rest of the paper the coefficient spaces will be
indexed by co-degree.  For this the notation $\Lambda(\delta,r,s)$ will be
used.  Here, $\delta$ stands for a flag that has been broken into $l$
segments by segmentating data $s$, and $r=(r_1,r_1,\ldots,r_l)$ provides a
\emph{co-degree} for each segment.  The same flag can perhaps be segmented
in many ways, evn if $l$ is fixed.  The boundary map preserves co-degree.

\section{Local-global homology}

The main definition of this paper can now be given.  Recall that if
$\delta$ is a flag on the convex polytope $\Delta$, and $s$ breaks
$\delta$ (and $V$) into $l$ segments, and if a multi-component co-degree
$r=(r_1,\dots,r_l)$ is given, then from all this a coefficient space
$\Lambda(\delta,r,s)$ has been defined.  Recall also that if $\delta'$ is
obtained from $\delta$ by removal of the $i$-th term from $\delta$, then
there is a natural map
\[
    \partial_i : \Lambda(\delta,r,s) \to \Lambda(\delta',r,s')
\]
such that the usual definition of $\bound$ will produce a boundary map on
\[
    A_r = (0 \to A_{r,n} \to A_{r,n-1} \to 
    \dots \to A_{r,1} \to A_{r,0} \to 0 )
\]
where $A_{r,i}$ is the direct sum of the $\Lambda(\delta_r)$ for
$\ord \delta=i$.

(The following detail is important.  When $\delta_i$ is removed from
$\delta$ to obtain $\delta'$, one might as a result have to change the
segmentation $s$.  This happens when $\delta_i$ is used to segment
$\delta$.  In this case the next term $\delta_{i+1}$ is used to segment
$\delta'$.  If $\delta_{i+1}$ is already used by $s$ to segment $\delta$,
or does not exist because $\delta_i$ is the last term in $\delta$, then
this component of the boundary map is treated as zero.  In this way, one
obtains either a satisfactory $s'$, or the zero map.)

Thus, for each co-degree $r=(r_1,\ldots,r_l)$, a complex $A_r$ has been
defined.  The boundary map $\bound$ of $A_r$ may cause the terms of the
segmenting flag $\delta_s$ to move rightwards.  This observation leads to
the following.  One can define subcomplexes of $A_r$ by placing conditions
of the segmentation subflags that are to be used.  These conditions must
of course allow the movement to the right of the segmentation terms.

One way to formulate this is to introduce a multi-dimension $s=(s_1 <
\dots < s_{l-1})$ and allow only those segmentations $\delta$ to be used,
for which the dimension $d=(d_1<\dots<d_l)$ is term by term at least as
large as $s$.  In this way one obtains for each $r$ and $s$ a complex
$A_{r,s}$.  Of course, if $s$ is too great compared to $r$, then the
complex will be zero.  The complex $A_{r,s}$ is the $A_w$ mentioned in
\S1.)

We now define the $(r,s)$ \emph{homology space} $H_{r,s}\Delta$ of
$\Delta$ to be the homology of $A_{r,s}$ at the level where $\ord\delta=i$
is equal to the total degree of $\Lambda(\delta,r,s)$.  In other words, it
is the homology at $A_{r,s;i}$, where $ i + r_1 + \dots + r_l = n$, for
$r$ gives the co-degree.  If $A_{r,s}$ exactly computes $H_{r,s}\Delta$
then its dimension $h_{r,s}\Delta$ is of course a linear function of the
flag vector.  The next section clarifies this.

As mentioned in the introduction, these spaces correspond to the
local-global intersection homology spaces of $\PDelta$, introduced in
\cite{bib.JF.LGIH}.

\section{Decorated bar diagrams}

The previous discussion made no use of coordinates, and does not give the
dimension of the various $\Lambda(\delta, r, s)$ spaces involved in the
construction of local-global homology.  This section provides another
approach.  The contribution made by a flag $\delta$ to the homology
$H_{r,s}\Delta$, and the maps between these contributions, can be
described using coordinates, via the use of decorated bar diagrams.  (Part
of the theory of bar diagrams was first published in the survey
paper~\cite{bib.MB.TASI}.)

Suppose, for example, that the dimension $n$ of $\Delta$ is eleven, and
that $\delta$ is a flag of dimension $d=(3 < 5 < 9 < 11)$. The
(undecorated) \emph{bar diagram}
\[
    \bardiagram{...|..|....|..}
\]
expresses this situation.  There are eleven dots, and a bar `\bardiagram{|}'
is placed after the $d_i$-th dots.  Each bar represents a term $\delta_i$ of
$\delta$.  Now choose a basis $e_1, \dots , e_n$ of $V$ such that for each
$i$, the initial sequence $e_1, \dots , e_j$ (with $j=d_i$) is a basis for
the span $\langle\delta_i\rangle$ of the vectors lying on the $i$-th face. 
Finally, let each dot represent not $e_i$ but the corresponding linear
function $\alpha_i$, that vanishes on each $e_j$ except $e_i$.  Call this
a \emph{system of coordinates} for $V$, that is \emph{subordinate} to
$\delta$.

By construction, each dot represent a linear function that vanishes on the
faces $\delta_i$ of $\delta$ represented by the bars that lie to its left. 
So that we can speak more concisely, we shall think of the dots and bars
as actually being the linear functions and terms of the flag respectively.

Each subset of the dots represents an element of the exterior algebra
generated by the linear functions on $V$.  To represent the selection of a
subset, promote the chosen dots `\bardiagram{.}' into circles `\bardiagram{o}'. 
Thus
\begin{equation}
\label{eqn:dbd} 
    \bardiagram{..o|o.|o.o.|oo}
\end{equation}
represents (or more concisely is) a degree six element of the exterior
algebra.  It is an example of a \emph{decorated bar diagram}.

(Each circle represents an element in the homology of the torus, that is
the generic or central fibre of the moment map.  As the fibre is moved
towards a face that lies to the left of the circle, so the circle shrinks
to a point.  Thus, circles represent $1$-cycles with specified vanishing
properties.  Similarly, the promotion of say $6$ dots to circles produces
a diagram that represent a $6$-cycle in the generic fibre, with specific
shrinking properties, as the fibre is moved to the faces of the flag. This
will be important, in the topological interpretation of exact homology.)

The exterior form (\ref{eqn:dbd}) has certain vanishing properties, with
respect to the faces of $\delta$.  The condition of \S4 is equivalent to
the following: That between each bar and the right hand end of the
diagram, there should be strictly more `\bardiagram{o}'s than `\bardiagram{.}'s. 
Our example satisfies this condition, and so is an \emph{admissable}
decorated bar diagram.

Now fix an unsegmented degree $r=(r_1)$ and define the coefficient space
$\Lambda(\delta,r)$ to be the span of the admissable $r$-circle bar
diagrams (considered as exterior forms).  (Because there is no
segmentation, $s$ is trivial, and will be omitted.)  For the
$\Lambda(\delta,r)$ to come together to produce a complex, the following
must hold.  First, $\Lambda(\delta,r)$ as a subspace of the degree~$r$
forms should not depend on the choice of a basis subordinate to $\delta$. 
Second, the boundary map should not depend on the basis (or in other words
should be covariant for such change). Third, when a bar is removed from an
admissable diagram, the result should also be admissable.  The last two
conditions are immediately seen to be true.

There are two ways to see that the first requirement (that
$\Lambda(\delta,r)$ not move when the subordinate basis is changed) is
true.  One method is to show that the span $\Lambda(\delta,s)$ of the
admissable diagrams is the solution set to a problem that can be
formulated without recourse to use of a basis.  This is the approach taken
in \S\S3--7.  The other method is to show directly that
$\Lambda(\delta,r)$ does not move, under change of subordinate basis.

Any change of subordinate basis can be obtained as a result of applying
the following moves.  First, one can multiply basis elements by non-zero
scalars.  Second, one can permute the basis elements (dots and circles),
provided so doing does not cause a dot or circle to pass over a bar. 
Thirdly, one can increase a basis element by some multiple of another
basis element, that lies to its right.

It is clear that applying either of the first two moves to the basis will
not change $\Lambda(\delta,r)$.  Regarding the third move, if a diagram
such as (\ref{eqn:dbd}) is admissable, then the result of moving one or
more `\bardiagram{o}'s to the right, perhaps over bars, will also be
admissable.  This is obvious, from the nature of the conditions defining
admissability.  The third type of move makes such changes.  Thus, the span
$\Lambda(\delta,r)$ of the admissable decorated diagrams does not move.

It has now been shown that the admissable bar diagrams, such as
(\ref{eqn:dbd}), define coefficient spaces $\Lambda(\delta,r)$ that can
be assembled to produce a complex.  It is left to the reader, to check that
it is exactly the same complex, as was defined in \S4.

The remainder of this section is devoted to the description of the
segmented form of the above construction.  As before, let
\[
    \bardiagram{...|..|....|..}
\]
denote a flag of dimension $d=(3 < 5 < 9 < 11)$, and now choose some of
the bars, say just the second, to segment the diagram.  The result
is
\[
    \bardiagram{...|..} \qquad \bardiagram{|....|..}
\]
or more concisely
\[
    \bardiagram{...|..!....|..}
\]
where the promotion of a `\bardiagram{|}' to a `\bardiagram{!}' indicates that is
is being used for segmentation.  Note that each segmenting face is also
used as the first face in the following segment.

As before, one can promote some of the `\bardiagram{.}'s to `\bardiagram{o}'s to
obtain a \emph{decorated bar diagram}, which will be \emph{admissable} if
between any bar (either `\bardiagram{|}' or `\bardiagram{!}') and the end of the
segment, there are more `\bardiagram{o}'s than `\bardiagram{.}'s.

To be able to assemble the span of the admissable diagrams into a complex,
certain requirements must be met.  They have already been formulated.  The
first is that the span should not depend on the choice of a subordinate
basis. Scalar and permutation moves on the basis clearly leave the span
unchanged, as in the single segment case.  Adding to one basis vector
another, lying to its right, has no effect on the span, provided the
second lies in the same segment as the first.

Now consider the situation, where one of the basis linear forms is
changed, by adding to it another basis linear form, as before lying to the
right, but this time in a different segment.  The way to have this have no
effect on the span is to have the dots and circles in a segmented diagram
represent not linear functions on $V$, but rather linear functions on the
(span of vectors lying on) the face at the right end of the segment (which
is $V$ for the last segment).

Thus, by having each decorated diagram (possibly segmented) represent an
element in a tensor product of exterior algebras, the span of all the
admissable diagrams becomes a space that does not move, when the
subordinate basis is changed.  This agrees with \S5.

For these spaces to be assembled into a complex, the boundary map must be
well defined.  As in the single segment situation, all is well when a
non-segmenting term (a `\bardiagram{|}' rather than a `\bardiagram{!}') is removed
from a decorated bar diagram.

The removal of a `\bardiagram{!}' terms has a more subtle effect, for the
segmentation will have to change.  Whatever rule is used, it must respect
admissability of decorated diagrams, and it must respect change of
subordinate basis.

The rule, which we state without prior justification, is this:  whenever
something such as
\[
    \bardiagram{!ooooo|}
\]
occurs in a decorated diagram, one can obtain a component of the boundary
by replacing it with
\[
    \bardiagram{ooooo!}
\]
while the result of removing any other `\bardiagram{!}' terms is zero.  
(The number of `\bardiagram{o}'s is not relevant, that one has
`\bardiagram{!}' followed by some circles, and then a `\bardiagram{|}'
is.)

It is clear that this rule will, from admissable diagrams, generate only
admissable diagrams.  This is because it can but only cause some extra
`\bardiagram{o}'s to appear, whenever the `\bardiagram{.}' and `\bardiagram{o}' counts
are to be compared.  The same is of course not true, for the replacement
of `\bardiagram{!.o.|}' by `\bardiagram{.o.!}' and similar situations.

The next task is to show that this part of the boundary respects change of
subordinate basis.  Given a fragment such as `\bardiagram{!.o.|}' that
contributes zero to this part of the boundary, whatever change of basis is
made, the contribution will still be zero.  Consider now a fragment such
as `\bardiagram{!ooo|}'. Scalar and permutation moves will have no effect
on the boundary.  This is obvious.  The only way that adding something
that lies to the right can have any effect is if that something lies to
the right of the whole fragment `\bardiagram{!ooo|}'.  This is because the
exterior algebra is antisymmetric.  Consider now the boundary contribution
`\bardiagram{ooo!}'.  As already described, the `\bardiagram{!}' induces a
restriction of linear functions, and so the just considered change of
basis will after all have no effect on the boundary map.

This concludes the presentation of the complex $A$ via decorated bar
diagrams.  The details of co-degree (number of `\bardiagram{.}'s in each
segment) and so forth are as in the earlier exposition~(\S5--7).  It is left
to the reader to verify that the two approaches lead to exactly the same
family $A_{r,s}$ of complexes.

\section{Summary and conclusions}

The flag vector of a convex polytope satisfies subtle linear inequalities,
and also non-linear inequalities, that are at present not known.  Exact
homology seems to be the only general method available to us, to prove
such results.  The difficulty with purely combinatorial means is firstly
that there are no natural maps, other than the deletion operators, between
the flags on a polytope, and secondly that the convexity of the polytope
has to be allowed to enter into the discussion in a significant way.  This
said, very few results relating to exact homology are known at present. 

For the usual middle perversity intersection homology theory, the three
approachs mentioned in the introduction are known to have significant
areas of agreement.  The recursive formula in \cite{bib.RS.GHV} for $h\PDelta$
can be unwound to express each Betti number as a sum of contributions due
to flags.  These contributions are exactly the same as those implicit in
the definitions of \S4.  The method of decorated bar diagrams provides a
reformulation of the recursive formula for $h\PDelta$, where now one
simply counts the number of valid diagrams.  As was indicated in \S8, once
the valid diagrams are known, it is not then hard to determine what the
corresponding space of coefficients should be.  Thus, the derivation from
\cite{bib.RS.GHV} of \S4 is without difficulty.  Finally, via the
identification of the exterior algebra $\bigwedge^\bullet V^*$ with the
homology of the generic fibre of the moment map, one can interpret the
complex in \S4 in terms of the construction of intersection homology
cycles on $\PDelta$, and relations between them.  This last will be
presented elsewhere.

The important property of intersection homology (for middle perversity
only) is that it seems to produce exact homology groups.  That there are
formula for such Betti numbers, and moreover as somewhat geometric
alternating sums, supports this view.  The local-global construction
defined in \S\S5--7 can, via the moment map, be translated into a class of
cycles on $\PDelta$, and relations between them.  These cycles and
relations have special properties with respect to the strata of $\PDelta$.
In this way one can translate the definition of local-global homology in
this paper into a topological one.  Indeed, if intersection homology were
not already known, it could have been discovered via its similarity to the
$H\Delta$ theory presented here.

There is a subtlety connected to the concept of a linear function of the
flag vector.  The flag vectors of convex polytopes span a proper subspace
of the space of all possible flag vectors.  Now, a linear function on a
subspace is not the same as a linear function on the whole space.  The
latter contains more information.  The constructions of this paper provide
linear functions of arbitrary (not necessarily polytope) flag vectors. 
As noted, they agree with the formulae (8)--(10) of~\cite{bib.JF.LGIH}.

The construction presented here of $H\Delta$ is probably not the only one.
In particular, it ought to be possible to define a theory, with the same
expected Betti numbers, but where $H\Delta$ is built out of the $H\delta$,
for all proper faces $\delta\subset\Delta$, together with perhaps a little
gluing information.  This can certainly be done in the simple case, where
it corresponds to the natural formula in that context for $h\Delta$ in
terms of the face vector.  For general polytopes, such a theory would
correspond to a linear function on arbitrary flag vectors, which agrees on
polytope flag vectors with the $h$-vector presented in this paper.  Such a
theory may be part of a geometric proof of exactness of homology.

One cannot, of course, prove that which is not true; and exact homology
cannot hold for a complex if the expected Betti number turns out to be
negative for some special polytope.  Bayer (personal communication) has an
example of a 5~dimensional polytope (the bipyramid on the cylinder on a
3-simplex) where this in fact happens.  This is discussed further in
\cite{bib.JF.LGIH}.  All this indicates that there are as yet unrealised
subtleties in the concept and proof of exact calculation.

The central definitions of this paper, namely of the coefficient spaces
$\Lambda(\delta,r,s)$ and the maps between them, use only linear algebra
and the deletion operator on flags.  Put another way, the largest part of
this paper has been the study of the flag
\[
    V_1 \subset V_2 \subset \dots \subset V_m \subset V
\]
of vector spaces associated to a flag $\delta$ of faces on $\Delta$, and
the vector spaces that can be defined from it.  Given the requirements of
exact homology, there is little else available that one could study. 
However, these definitions are very closely linked to the usual
intersection homology theory and also its local-global variant.  These
connections indicate that there may be unremarked subtleties in the linear
algebra of a flag, and undiscovered simplicity in intersection homology. A
better understanding of the intersection homology of Schubert varieties
would be very useful.

Finally, as mentioned in the Introduction, one could wish for a similar
theory that applies to $i$-graphs, and similar combinatorial objects. 
This will probably require a different sort of linear algebra.


\begin{thebibliography}{99}

\bibitem{bib.MB.TASI}
M.M. Bayer, Face numbers and subdivisions of convex polytopes, in
{\it POLYTOPES: Abstract, Convex and Computational} (ed.\
T.~Bisztriczky et al), (1994) 155--171. Kluwer.

\bibitem{bib.PD.WC1}
P. Deligne, La conjecture de Weil, I, {\it Publ. Math. IHES} {\bf 43}
(1974), 273--307

\bibitem{bib.PD.WC2}
\bibrule, La conjecture de Weil, II, {\it Publ. Math. IHES} {\bf 52}
(1980), 137--252

\bibitem{bib.JD-FL.IHNP}
J. Denef and F. Loeser, Weights of exponential sums, intersection
homology and Newton polyhedra, {\it Invent. Math.} {\bf 106} (1991),
275--294

\bibitem{bib.KF.IHTV}
K. Fieseler, Rational intersection homology of projective toric
varieties, {\it J. reine angew. Math.} {\bf 413} (1991), 88--98

\bibitem{bib.JF.QTHGFV}
J. Fine, Quantum topology, hypergraphs and flag vectors,
preprint q-alg/9708001 (August 1997)

\bibitem{bib.JF.SFV}
\bibrule, Shelling and flag vectors,
preprint q-alg/9710002 (October 1997)

\bibitem{bib.JF.LGIH}
\bibrule, Local-global intersection homology, 
preprint alg-geom/9709011 (September 1997)

\bibitem{bib.RS.GHV}
R. Stanley,
Generalized $h$-vectors, intersection cohomology of
toric varieties and related results, {\it Adv. Stud. Pure Math.} {\bf
11} (1987), 187--213

\end{thebibliography}
\end{document}